\documentclass[10pt,reqno]{amsart}
\usepackage{eucal,fullpage,times,amsmath,amsthm,amssymb,mathrsfs,stmaryrd,color,enumerate,accents}
\usepackage[all]{xy}
\usepackage{url}

\newcommand{\calX}{{\mathcal{X}}}

\newcommand{\calO}{{\mathcal{O}}}

\newcommand{\calA}{{\mathcal{A}}}

\newcommand{\calF}{\mathcal{F}}

\newcommand{\calC}{\mathcal{C}}

\newcommand{\A}{\mathbf{A}}

\newcommand{\Z}{\mathbf{Z}}
\newcommand{\N}{\mathbf{N}}
\newcommand{\C}{\mathbf{C}}
\newcommand{\F}{\mathbf{F}}
\newcommand{\Q}{\mathbf{Q}}
\renewcommand{\P}{\mathbf{P}}

\newcommand{\Spec}{{\mathrm{Spec}}}

\newcommand{\Hom}{\mathrm{Hom}}
\newcommand{\Fun}{\mathrm{Fun}}

\newcommand{\Shv}{\mathrm{Shv}}

\newcommand{\Tor}{\mathrm{Tor}}

\newcommand{\Sym}{\mathrm{Sym}}

\newcommand{\Mod}{\mathrm{Mod}}

\newcommand{\Alg}{\mathrm{Alg}}

\newcommand{\Set}{\mathrm{Set}}

\newcommand{\R}{\mathrm{R}}

\newcommand{\Ch}{\mathrm{Ch}}

\newcommand{\zar}{\mathrm{Zar}}
\newcommand{\an}{\mathrm{an}}

\newcommand{\im}{\mathrm{im}}
\newcommand{\id}{\mathrm{id}}
\newcommand{\ev}{\mathrm{ev}}

\newcommand{\Fil}{\mathrm{Fil}}

\newcommand{\Bl}{\mathrm{Bl}}

\newcommand{\gr}{\mathrm{gr}}
\newcommand{\dR}{\mathrm{dR}}

\renewcommand{\inf}{\mathrm{inf}}

\newcommand{\opp}{\mathrm{op}}

\renewcommand{\ker}{\mathrm{ker}}

\newcommand{\Cech}{\mathrm{Cech}}
\newcommand{\Tot}{\mathop{\mathrm{Tot}}}
\newcommand{\Comp}{\mathrm{Comp}}

\newcommand{\Ho}{\mathrm{Ho}}

\newcommand{\comment}[1]{}

\newcommand{\cosimp}[3]{\xymatrix@1{#1 \ar@<.4ex>[r] \ar@<-.4ex>[r] & {\ }#2 \ar@<0.8ex>[r] \ar[r] \ar@<-.8ex>[r] & {\ } #3 \ar@<1.2ex>[r] \ar@<.4ex>[r] \ar@<-.4ex>[r] \ar@<-1.2ex>[r] & \cdots }}

\begin{document}
\bibliographystyle{alpha}
\newtheorem{theorem}{Theorem}[section]
\newtheorem*{theorem*}{Theorem}
\newtheorem*{condition*}{Condition}
\newtheorem*{definition*}{Definition}
\newtheorem*{corollary*}{Corollary}
\newtheorem{proposition}[theorem]{Proposition}
\newtheorem{lemma}[theorem]{Lemma}
\newtheorem{corollary}[theorem]{Corollary}
\newtheorem{claim}[theorem]{Claim}

\theoremstyle{definition}
\newtheorem{definition}[theorem]{Definition}
\newtheorem{question}[theorem]{Question}
\newtheorem{remark}[theorem]{Remark}
\newtheorem{guess}[theorem]{Guess}
\newtheorem{example}[theorem]{Example}
\newtheorem{condition}[theorem]{Condition}
\newtheorem{warning}[theorem]{Warning}
\newtheorem{notation}[theorem]{Notation}
\newtheorem{construction}[theorem]{Construction}

\title{Completions and derived de Rham cohomology}
\begin{abstract}
We show that Illusie's derived de Rham cohomology (Hodge-completed) coincides with Hartshorne's algebraic de Rham cohomology for a finite type map of noetherian schemes in characteristic $0$; the case of lci morphisms was a result of Illusie. In particular, the $E_1$-differentials in the derived Hodge-to-de Rham spectral sequence for singular varieties are often non-zero. Another consequence is a completely elementary description of Hartshorne's algebraic de Rham cohomology: it is computed by the completed Amitsur complex for any variety in characteristic $0$.
\end{abstract}
\author{Bhargav Bhatt}
\maketitle

\section{Introduction}
\label{sec:intro} 

Let $X$ be an algebraic variety over a characteristic $0$ field $k$.  If $X$ is smooth, then the cohomology of the algebraic de Rham complex $\Omega^*_{X/k}$ is a good cohomology theory for $X$ \cite{GrothendieckAlgdR}: it is canonically isomorphic to the Betti cohomology of $X^\an$ if $k = \C$, and hence gives a purely algebraic definition of the ``true'' cohomology groups of $X$.  However, if $X$ is singular, then the de Rham complex of $X$ is rather poorly behaved as $\Omega^1_{X/k}$ is not locally free, and the cohomology of $\Omega^*_{X/k}$ does not compute the desired groups (see \cite[Example 4.4]{ArapuraKangdR}).  Instead, motivated by the analogy between formal completions in algebraic geometry and tubular neighbourhoods in topology, Hartshorne \cite{HartshorneAlgdR} discovered\footnote{This construction is due independently to Deligne (unpublished, 1970), and Herrera-Lieberman \cite{HerreraLieberman}.} the ``correct'' recipe: if $i:X \hookrightarrow Y$ is a closed immersion with $Y$ smooth,  then the cohomology of the formal completion of the de Rham complex of $Y$ along $X$ is independent of the choice of $Y$, and computes the ``true'' cohomology groups of $X$. We call this theory ``algebraic de Rham cohomology'' in this note.

Another generalisation of de Rham theory for singular varieties, suggested by Illusie \cite[\S VIII]{IllusieCC2}, entails replacing the cotangent sheaf $\Omega^1_{X/k}$ with the cotangent complex $L_{X/k}$ in the definition of the de Rham complex; we refer to the resulting cohomology theory as ``(Hodge-completed) derived de Rham cohomology.'' Essentially by direct calculation of the cotangent complex, Illusie showed \cite[Corollary VIII.2.2.8]{IllusieCC2} that derived de Rham cohomology coincides with algebraic de Rham cohomology for varieties with lci singularities. Our main goal in this note is to show that this last comparison holds in complete generality:

\begin{theorem*}
The Hodge-completed derived de Rham cohomology of any finite type morphism of noetherian $\Q$-schemes is canonically isomorphic to Hartshorne's algebraic de Rham cohomology (ignoring filtrations).
\end{theorem*}

For a precise formulation, see Corollary \ref{cor:ddrgeomcomp}. One consequence is a manifestly choice-free construction of algebraic de Rham cohomology where certain properties (like the Kunneth formula) are obvious. Moreover, the Hodge filtration on derived de Rham cohomology defines a new filtration --- the derived Hodge filtration --- on algebraic de Rham cohomology. This filtration is finer than the standard Hodge-theoretic filtrations on algebraic de Rham cohomology (such as the infinitesimal Hodge filtration or the Hodge-Deligne filtration), and is a shadow of a spectral sequence --- the derived Hodge-to-de-Rham spectral sequence --- computing algebraic de Rham cohomology. The theorem above shows the {\em non-vanishing} of many differentials of this spectral sequence: the cohomology of the graded pieces of the derived de Rham complex is often infinite dimensional (due to singularities), but the cohomology of the total complex is always finite by the above theorem, so there must be ``cancellation'' throughout the spectral sequence. 

An elementary consequence of our {\em methods} is an explicit description of Betti cohomology for algebraic varieties via the completed Amitsur complex, generalising Grothendieck's recipe \cite{GrothendieckCrystals} from the smooth case. We present it here in the absolute affine setting, and simply note that it naturally extends to arbitrary flat maps (see Corollary \ref{cor:elemdrglob}):

\begin{corollary*}
Let $X = \Spec(R)$ be a finite type affine $\C$-scheme. Then there is a multiplicative isomorphism
\begin{equation}
	\label{eq:formalperiod}
	\R\Gamma(X^\an,\C) \simeq \Big( R \to \widehat{R^{\otimes 2}} \to \widehat{R^{\otimes 3}} \to \dots \Big)
\end{equation}
where the completion on the right takes place along the multiplication maps $R^{\otimes n} \to R$, and the differentials are the alternating sums of the coprojection maps. 
\end{corollary*}

Isomorphism \eqref{eq:formalperiod} above can be viewed as a ``formal'' period integral. For example, if $R = \C[t,t^{-1}]$ is the co-ordinate ring of a torus $T$, then the closed form $\frac{dt}{t} \in \Omega^1_{T/\C}$ is a generator of $H^1(T^\an,\C)$ under the de Rham-Betti isomorphism, and is represented by the cocycle 
\[ \int_{t \otimes 1}^{1 \otimes t} \frac{dt}{t}:= \log(\frac{1 \otimes t}{t \otimes 1}) \in \widehat{R \otimes_\C R}\] 
under isomorphism \eqref{eq:formalperiod} above.

Finally, we remark that the above theorem does not admit an evident positive characteristic analogue outside the lci setting: derived de Rham cohomology (completed or not) of a singular algebraic variety in characteristic $p$ is often non-zero in arbitrarily {\em negative} cohomological degrees (see \cite[Example 3.21]{Bhattpadicddr}), while classical cohomology theories live only in non-negative degrees. However,  the corollary above does have a meaningful extension (see Remark \ref{rmk:poscharamitsur}).

\subsection*{A sketch of the proof of the theorem} Since both derived de Rham cohomology and algebraic de Rham cohomology are well-behaved (and isomorphic) for smooth morphisms, the key case is that of closed immersions. The local statement is: if $I \subset A$ is an ideal in a noetherian $\Q$-algebra $A$, then the $I$-adic completion $\widehat{A}$ coincides with the (Hodge-completed) derived de Rham cohomology $\widehat{\dR}_{(A/I)/A}$ of $A \to A/I$. When $I$ is cut out by a regular sequence, this is essentially immediate: the Hodge-truncations of $\widehat{\dR}_{(A/I)/A}$ are naturally identified with the $I$-adic truncations of $A$ lifting the canonical isomorphism $L_{(A/I)/A}[-1] \simeq I/I^2$. However, if $I$ is not regular, then these two filtrations have vastly different graded pieces (as the cotangent complex is often left-unbounded), so we must look elsewhere for a comparison. We solve this problem by introducing a third player: the Adams completion $\Comp_A(A,I)$ of $A$ along $I$. The construction of $\Comp_A(A,I)$  is borrowed from the Adams resolution in algebraic topology and is due to G. Carlsson. Using a convergence theorem of Quillen, we endow $\Comp_A(A,I)$ with a complete filtration whose truncations coincide with the Hodge-truncations of $\widehat{\dR}_{(A/I)/A}$; taking limits then yields a filtered isomorphism $\Comp_A(A,I) \simeq \widehat{\dR}_{(A/I)/A}$. On the other hand, a theorem of Carlsson also identifies $\Comp_A(A,I)$ with the classical $I$-adic completion $\widehat{A}$ (ignoring filtrations) for noetherian rings, thereby proving the claim.

\subsection*{Organisation of this note} Homological conventions (especially surrounding homotopy-limits) are discussed in the next paragraph. In \S \ref{sec:cosimpalg}, we formulate and prove a descent statement for the cotangent complex. The aforementioned Adams completion operation is studied in \S \ref{sec:adamscomp}. The main comparison theorems are proven in \S \ref{sec:mainthm}, and the resulting web of filtrations (together with their relation to algebraic cycles) is discussed in \S \ref{sec:derhodgefilt}. 
%The last section \S \ref{sec:luriecomp} studies other completion operations for rings and compares them to the Adams completion; this can be read following \S \ref{sec:adamscomp} and is independent of the rest of the paper.

\subsection*{Notation and conventions}
For each $n \in \Z_{\geq 0}$, let $[n]$ denote the poset $\{ 0 < 1 < \dots < n \}$, and let $\Delta$ be the simplex category. For any category $\calC$, we use $s\calC$ and $c\calC$ to denote the categories of simplicial and cosimplicial objects respectively; for both cosimplicial and simplicial objects $A$ over $\calC$, we write $A_n := A([n])$ unless otherwise specified. If $\calC$ has a final object, we use $\calC_\ast$ to denote the category of pointed objects in $\calC$, i.e., the category of maps $\ast \to X$ in $\calC$ with $\ast$ a fixed final object; if the initial and final objects in $\calC$ coincide, then $\calC_* \simeq \calC$, and $\calC$ is called {\em pointed}.  For a pointed category $\calC$, a map $f:X \to Y$ of objects in $c\calC$ is called {\em null-homotopic} if it is homotopic to the constant map (defined via the final object); an object $X \in c\calC$ is said to be {\em homotopy-equivalent to $0$} if $\id_X$ is null-homotopic. 

We use model categories to discuss homotopy-limits. Given a (nice) model category $\calC$ and a small indexing category $I$, we write $R\lim_I:\Fun(I,\calC) \to \calC$ for the homotopy-limit functor, which is a right Quillen adjoint to the ``constant'' functor $\calC \to \Fun(I,\calC)$ for the injective model structure on the latter; alternatively, $R\lim_I$ is the derived pushforward on $\calC$-valued sheaves along the map of sites $I \to \ast$ (with the indiscrete topology). We will implicitly use the commutativity of homotopy-limits with other homotopy-limits, i.e., the identification of the composites $\Fun(I \times J,\calC) \stackrel{R\lim_J}{\to} \Fun(I,\calC) \stackrel{R\lim_I}{\to} \calC$ and $\Fun(I \times J,\calC) \stackrel{R\lim_I}{\to} \Fun(J,\calC) \stackrel{R\lim_J}{\to} \calC$ for small indexing categories $I$ and $J$. The relevant examples are: $I = \N^\opp$ (see below), $J = \Delta$, and $\calC = \Ch(\calA)$ for some Grothendieck abelian category $\calA$ (via the ``Spaltenstein'' model structure, see \cite[Proposition 1.3.5.3]{LurieHA}). We follow the convention $\Tot := R\lim_\Delta$; note that if $K$ is a cosimplicial abelian group, then $\Tot(K)$ is equivalent to the chain complex underlying $K$ via the Dold-Kan correspondence.

Given a category $\calC$ with finite colimits and an object $A \in \calC$, we write $\Cech(A \to \cdot):\calC_{A/} \to c\calC_{A/}$ for the ``Cech conerve'' functor, i.e., the left adjoint to the forgetful functor $\ev_{[0]}:c\calC_{A/} \to \calC_{A/}$. Explicitly, for any map $f:A \to B$, one has $\Cech(A \to B)([m]) \simeq \sqcup_{[m]} B$ where the coproduct is computed relative to $A$. If $\calC$ is also a model category, then the functor $\Cech(A \to \cdot)$ is a left Quillen adjoint under the injective model structure, and we abusively write use $\Cech(A \to \cdot)$ to denote the left derived functor as well. Explicitly, $\Cech(A \to B)([m])$ is computed as $\sqcup_{[m]} P$ where $P$ is a cofibrant replacement for $B$ in $c\calC_{A/}$, and the coproduct is computed relative to $A$.

We write $\N^\opp$ for the category underlying the poset underlying the non-positive natural numbers under the usual ordering. Given a diagram $K:\N^\opp \to \calC$ in some category $\calC$, we use $\{K_k\}$ to denote the corresponding object of $\Fun(\N^\opp,\calC)$, and similarly for objects in $D(\Fun(\N^\opp,\calA))$ when $\calA$ is abelian. For a Grothendieck abelian category $\calA$, we informally refer to an object $\{K_k\} \in D(\Fun(\N^\opp,\calA))$ as the associated system of quotients for an $\N^\opp$-indexed complete separated filtration on $\widehat{K} := R\lim K_k \in D(\calA)$; this convention may be formalised by viewing the auxilliary system $\{\ker(\widehat{K} \to K_k)\} \in D(\Fun(\N^\opp,\calA))$ of homotopy kernels as a filtration on $\widehat{K}$, but we never do this. A cochain complex $K$ over $\calA$ defines an object $\{K_n\} \in D(\Fun(\N^\opp,\calA))$ via $K_n :=  K/\sigma^{\geq n} K$, where $\sigma^{\geq n}$ denotes the stupid truncation of $K$ in cohomological degrees $\geq n$; one has $K \simeq R\lim_n K_n$, and we refer to the resulting filtration on $K$ as the stupid filtration on $K$.

The category $s\Alg$ (resp. $cs\Alg$) of simplicial (resp. cosimplicial simplicial) commutative rings plays a key role in this note. A map in $s\Alg$ is called {\em surjective} if it is surjective on $\pi_0$. We implicitly use the following: if $A \to B$ is an equivalence in $s\Alg$ and $M$ is a termwise flat simplicial $A$-module, then $M \simeq M \otimes_A B$. We identify ordinary rings with the corresponding constant simplicial or cosimplicial rings. For $A \in s\Alg$, we write $s\Mod_A$ for the category of simplicial $A$-modules. For $A \in cs\Alg$, we write $cs\Mod_A$ or $c(s\Mod_A)$ for the category cosimplicial simplicial $A$-modules $M$, i.e., the data of a simplicial $A_n$-module $M_n$ for each $n \in \Delta$ together with cosimplicial structure maps $M_n \to M_m$ linear over $A_n \to A_m$ for each map $[n] \to [m]$ in $\Delta$. There are well-defined interval objects $\Delta[1] \in c(s\Mod_A)$ (induced by the ones in $c\Mod_\Z$ by base change), so it makes sense say that an object $K \in c(s\Mod_A)$ is homotopy-equivalent to $0$. A chain complex $K$ over some abelian category is called {\em connective} if $\pi_i(K) = H^{-i}(K) = 0$ for $i < 0$, and {\em coconnective} if $\pi_i(K) = H^{-i}(K) = 0$ for $i > 0$.

\subsection*{Acknowledgements} I thank Davesh Maulik for numerous useful conversations, one of which led to the writing of this note. In addition, I am grateful to Sasha Beilinson, Johan de Jong,  Luc Illusie, and Andrew Snowden for enlightening comments, questions, and encouragement. Thanks are also due to Srikanth Iyengar for making me aware of the origins of the Greenlees-May completion.

%Thanks are also due to Anatoly Preygel for a discussion leading to Example \ref{ex:completionzpbar}, and to Srikanth Iyengar for making me aware of the origins of the Greenlees-May completion from \S \ref{sec:luriecomp}.

\section{Some cosimplicial algebra}
\label{sec:cosimpalg}

Our goal in this section is to prove Corollary \ref{cor:cotcomplexcosimpalg}, a slightly more general version of \cite[Example 2.16]{BhattdeJong}.

\begin{construction}
\label{cons:cosimpobj}
Let $\calC$ be a category with finite coproducts. For each $X \in \calC$, define $\Delta_X(\bullet) := \Cech(\emptyset \to X) \in c\calC$. Explicitly, one has
\[ \Delta_X([n]) = \coprod_{i \in [n]} X \]
for each $[n] \in \Delta$. More canonically, there is a functor $\Set^f \to \calC$ defined via $S \mapsto \coprod_{s \in S} X$, and the cosimplicial object $\Delta_X(\bullet)$ is the image of the standard cosimplicial set $\Delta(\bullet) \in c\Set^f$ under the induced functor.
\end{construction}

\begin{proposition}
\label{prop:cosimpobjnull}
Let $\calC$ be a pointed category closed under finite colimits. For any $X \in \calC$, the object $\Delta_X(\bullet) \in c\calC$ is homotopy-equivalent to $0$.
\end{proposition}
\begin{proof}
Consider the functors
\[ \Set^f \stackrel{F}{\to} \Set^f_* \stackrel{G}{\to} \calC\]
defined via $F(S) = (S \coprod \{\ast\},\ast)$, and 
\[ G((T,t_0)) = \big(\coprod_{t \in T} X\big) / X\]
where $X \to \coprod_{t \in T} X$ is the map coming from the base point $t_0 \in T$. The object $\Delta_X(\bullet) \in c\calC$ is obtained by applying $(G \circ F)$ to the usual cosimplicial set $\Delta(\bullet) \in c\Set^f$. Hence, it suffices to check that $F(\Delta(\bullet)) \in c\Set^f_*$ is homotopy-equivalent to $0$. For this, recall that $\Delta(\bullet) = \Cech(\emptyset \to \{1\}) \in c\Set^f$, i.e., $\Delta([n]) = \sqcup_{i \in [n]} \{1\} = [n]$.  Application of $F$ commutes with formation of conerves (as $F$ is left adjoint to the forgetful functor $\Set^f_* \to \Set^f$), so $F(\Delta(\bullet)) = \Cech(f)$ where $f:(\{\ast\},\ast) \to (\{1,\ast\},\ast))$ is the unique map. Since $f$ has a section in $\Set^f_*$, it is standard to check that $F(\Delta(\bullet))$ is homotopy-equivalent to $0$ in $c\Set^f_*$. 
\end{proof}

\begin{construction}
\label{cons:cosimpmod}
Consider the cosimplicial abelian group $\Delta_\Z(\bullet) \in c\Mod_{\Z} \subset c(s\Mod_{\Z})$ defined via Construction \ref{cons:cosimpobj} applied to $\calC = \Mod_{\Z}$. For any cosimplicial simplicial ring $A$, we define $\Delta_A(\bullet) \in c(s\Mod_A)$ via $\Delta_{\Z}(\bullet) \otimes_{\Z} A$, where the base change takes place along the natural map $\Z \to A$. Explicitly, $\Delta_A(\bullet)$ is simply $\oplus_{i \in [n]} A_n$ in cosimplicial degree $n$; if each $A_n$ is an ordinary ring, then $\Delta_A(\bullet) \in c\Mod_A \subset cs\Mod_A$.
\end{construction}

\begin{corollary}
\label{prop:cosimpmodnull}
Let $A$ be a cosimplicial simplicial ring. Then $\Delta_A(\bullet)$ is homotopy equivalent to $0$ in $c(s\Mod_A)$.
\end{corollary}
\begin{proof}
The case $A = \Z$ follows from Proposition \ref{prop:cosimpobjnull}, and the general case follows by base change.
\end{proof}

\begin{lemma}
\label{lem:hoteq0}
Let $A$ be a cosimplicial simplicial ring. The forgetful functor $c(s\Mod_A) \to s\Mod_{A_0}$ has a left adjoint $F$. For any $M \in s\Mod_{A_0}$,  the object $F(M)$ is homotopy-equivalent to $0$ in $c(s\Mod_A)$.
\end{lemma}
\begin{proof}
The existence of the adjoint follows from the adjoint functor theorem. Explicitly, the formula is
\[ F(M)_n = \bigoplus_{i \in \Hom([0],[n])} A(i)^* M \simeq \bigoplus_{i \in [n]} A(i)^* M, \]
where $A(i):A_0 \to A_n$ is the cosimplicial structure map defined by $i$, and the cosimplicial structure of $F(M)$ is the evident one. An explicit homotopy between the identity map on $F(M)$ and the $0$ map can be constructed using the formulas in \cite[Example 2.16]{BhattdeJong}. Alternatively, observe that the formula for $F$ can be rewritten as
\[ F(M)_n = \bigoplus_{i \in [n]} A(i)^* M = \bigoplus_{i \in [n]} M \otimes_{A_0} A_0 \otimes_{A_0,A(i)} A_n = M \otimes_{A_0} \Big(\bigoplus_{i \in [n]} A_0 \otimes_{A_0,A(i)} A_n \Big) = M \otimes_{A_0} F(A_0). \]
In other words, $F(A_0)$ defines a $A_0$-module in the category $cs\Mod_A$, and $F$ is the unique colimit preserving functor prescribed by setting $F(A_0)$ as above. Hence, it suffices to check that $F(A_0)$ is homotopy-equivalent to $0$. However, $F(A_0)$ is easily identified with $\Delta_{A}(\bullet)$, so the claim follows from Corollary \ref{prop:cosimpmodnull}.
\end{proof}

The next lemma describes the effect of taking (derived) iterated tensor products on differential forms; it is formulated in terms of derived Cech conerves of maps in $s\Alg$, a notion discussed further in Construction \ref{cons:adamscomp}.

\begin{lemma}
\label{lem:cotcomplexcosimpalg}
Let $A \in s\Alg$, and let $P$ be a simplicial polynomial $A$-algebra. Let $P(\ast) = \Cech(A \to P) \in c(s\Alg_{A/})$ be the Cech conerve. Then $\Omega^k_{P(\ast)/A} \in c(s\Mod_{P(\ast)})$ is homotopy-equivalent to $0$ for $k > 0$.
\end{lemma}
\begin{proof}
It suffices to show the $k = 1$ case as homotopy-equivalences are preserved under the termwise application of the exterior power functor. The Kunneth formula gives $\Omega^1_{P(\ast)/A} \simeq F(\Omega^1_{P/A})$ where $F:s\Mod_P \to c(s\Mod_{P(\ast)})$ is the left adjoint to the restriction $c(s\Mod_{P(\ast)}) \to s\Mod_P$, so the claim follows from Lemma \ref{lem:hoteq0}.
\end{proof}

The following observation is the main result of this section. 

\begin{corollary}
\label{cor:cotcomplexcosimpalg}
Fix a map $A \to B$ in $s\Alg$. Then
\[ \Tot(\wedge^k L_{\Cech(A \to B)/A}) \simeq 0 \quad \mathrm{and} \quad \Tot(\wedge^k L_{B/\Cech(A \to B)}) \simeq \wedge^k L_{B/A}.\]
\end{corollary}
\begin{proof}
Choose a simplicial polynomial $A$-algebra resolution $P \to B$ with Cech conerve $P(\ast) \in c(s\Alg_{A/})$, and let $I \subset P(\ast)$ is the cosimplicial simplicial ideal defining the augmentation $P(\ast) \to P$. Then we have the following explicit models:
\[ \wedge^k L_{\Cech(A \to B)/A} = \Omega^k_{P(\ast)/A} \in c(s\Mod_{P(\ast)}) \quad \mathrm{and} \quad  \wedge^k L_{B/\Cech(A \to B)} =  \wedge^k (I/I^2[1])   \in c(s\Mod_{P}).\]
The first claim then follows from Lemma \ref{lem:cotcomplexcosimpalg}. The same lemma also implies that $\Omega^1_{P(\ast)/A} \otimes_{P(\ast)} P$ and its wedge powers are homotopy-equivalent to $0$ in $c(s\Mod_P)$. The second claim now follows by considering $\wedge^k$ of the (cosimplicial) transitivity triangle
\[ \Omega^1_{P(\ast)/A} \otimes_{P(\ast)} P \to \Omega^1_{P/A} \to I/I^2[1] \]
of simplicial $P$-modules.
\end{proof}

\begin{remark}
Corollary \ref{cor:cotcomplexcosimpalg} yields the following (derived) descent statement: for a fixed base ring $k$, the functor $L_{-/k}$ is a sheaf for the flat topology on $\Alg_{k/}$. In fact, Corollary \ref{cor:cotcomplexcosimpalg} shows that 
\[ \Tot(L_{A/k} \otimes_A \Cech(A \to B)) \simeq \Tot(L_{\Cech(A \to B)/k})\]
for any map $f:A \to B$. If $f$ is faithfully flat, then flat descent \cite{GrothendieckFlatDescent} identifies the left hand side with $L_{A/k}$, which proves the descent assertion above. Indeed, if $M$ is a discrete $A$-module (such as $\pi_i(L_{A/k})$), then $\Tot(M \otimes_A \Cech(A \to B)) \simeq M$ by Grothendieck's argument for flat descent of {\em modules}; the claim follows from this one by considering the Postnikov filtration on the cosimplicial complex $L_{A/k} \otimes_A \Cech(A \to B)$. In a sequel, we will use this observation to define and study the derived de Rham cohomology of an Artin stack (in any characteristic).
\end{remark}

\section{The Adams completion}
\label{sec:adamscomp}

\subsection{Generalities}

Our goal in this section is to introduce a variant of the usual $I$-adic completion operation on commutative rings $A$ with ideals $I$. More precisely, we study the following construction of derived Cech conerves (whose name is inspired by the Adams resolution in algebraic topology):

\begin{construction}
\label{cons:adamscomp}
Let $f:A \to B$ be a map in $s\Alg$. Then we use $\Cech(A \to B) \in cs\Alg$ to denote a (derived) Cech conerve of $f$; this is computed by the usual Cech conerve of the map $A \to P$ for any simplicial polynomial $A$-algebra resolution $P \to B$, and is independent up to homotopy of this choice. For any $M \in \Ch(A)$, we define the {\em Adams completion} of $M$ along $f$ as 
\[ \Comp_A(M,f) = \Tot(M \otimes_A \Cech(A \to B)) \in \Ch(A).\] 
If $M = C$ for some $C \in s\Alg_{A/}$, then $\Comp_A(C,f) \simeq \Comp_A(C,f \otimes_A \id_C)$ is naturally an $E_\infty$-algebra. If $f$ is a quotient map of ordinary rings with kernel $I \subset A$, then we also write $\Comp_A(M,I)$ for $\Comp_A(M,A \to A/I)$.
\end{construction}

\begin{remark}
In Construction \ref{cons:adamscomp}, we may use any simplicial flat $A$-algebra resolution $P \to B$ to compute $\Comp_A(M,f)$. Moreover, to get a completely functorial construction, one may use the canonical resolution.
\end{remark}

\begin{remark}
Assume $f$ is surjective in the setup of Construction \ref{cons:adamscomp}. By \cite[Corollary 6.7]{CarlssonDercomp}, for any $M \in s\Mod_A \subset \Ch(A)$, the totalisation of $M \otimes_A \Cech(A \to B) \in cs\Mod_A$ computed either in $s\Mod_A$ or in $\Ch(A)$ is the same. In particular, $\Comp_A(A,f)$ is naturally an object of $s\Alg_{A/}$.
\end{remark}

\begin{example}
\label{ex:cechnervesection}
Say $f:A \to B$ is a map in $s\Alg$ with a section. Then $\Comp_A(M,f) \simeq M$ for any $M \in \Ch(A)$. In fact, the section gives a homotopy-equivalences $A \simeq \Cech(A \to B)$  and $M \simeq  M \otimes_A \Cech(A \to B)$ in $c\Ch(A)$.
\end{example}

\begin{example}
\label{ex:cechnervefppf}
Say $f:A \to B$ is a faithfully flat map in $\Alg \subset s\Alg$. Then $\Comp_A(M,f) \simeq M$ for any $M \in \Ch(A)$. Indeed, this is Grothendieck's argument \cite{GrothendieckFlatDescent} for flat descent: this follows from Example \ref{ex:cechnervesection} when $f$ has a section, and the general case is reduced to this one by base changing along the faithfully flat map $A \to B$.
\end{example}

\begin{example}
\label{ex:cechnerveopen}
Let $f:A \to B$ be a (derived) epimorphism in $s\Alg$,  i.e., assume that multiplication induces an equivalence $B \otimes_A B \simeq B$; examples include localisations as well as some large quotients, such as $\C[\Q_{\geq 0}] \to \C$ or $\overline{\Z_p} \to \overline{\F_p}$. Then $\Cech(A \to B) \simeq B$, and hence $\Comp_A(M,f) \simeq M \otimes_A B$ for any $M \in s\Mod_A$.
\end{example}

\subsection{Identifying the Adams completion of a noetherian ring}

Following Carlsson, we will show that the Adams completion of a quotient map $A \to A/I$ coincides with the $I$-adic completion of $A$ when $A$ is noetherian.

\begin{definition}
Fix a Grothendieck abelian category $\calA$. An object $\{M_n \} \in \Fun(\N^\opp,\calA)$ is {\em strict-essentially $0$} if there exists an integer $k > 0$ such that the maps $M_n \to M_m$ are $0$ for all $n - m \geq k$. An object $\{K_n\} \in D(\Fun(\N^\opp,\calA))$ is {\em strict-essentially $0$} if the system $\{H^i(K_n)\} \in \Fun(\N^\opp,\calA)$ is strict-essentially $0$ for each $i$. 
\end{definition}

\begin{remark}
Fix a Grothendieck abelian category $\calA$. If $\{K_n\} \in D(\Fun(\N^\opp,\calA))$ is strict-essentially $0$, then $R\lim_n K_n \simeq 0$. In fact, the full subcategory of $D(\Fun(\N^\opp,\calA))$ spanned by strict-essentially $0$ objects is a thick triangulated subcategory, and the $R\lim_n$-functor factors through the quotient category (often called the derived category of pro-systems in $\calA$).
\end{remark}

\begin{example}
	\label{ex:artinrees}
Let $N \subset M$ be an inclusion of finitely generated modules over a noetherian ring $R$, and let $I \subset R$ be an ideal. Consider the map  $f:\{N/I^nN\} \to \{N/(I^nM \cap N)\}$ in $\Fun(\N^\opp,\Mod_R)$. Then $f$ is surjective, and the kernel $\{(I^nM \cap N)/ I^nN\}$ is strict-essentially $0$ by the classical Artin-Rees lemma.
\end{example}

We show that strict-essentially $0$ systems form an ideal:

\begin{lemma}
\label{lem:sezerotensorproduct}
Let $A$ be a ring. Let $\{K_n\} \in D^{\leq 0}(\Fun(\N^\opp,\Mod_A))$ be a strict-essentially $0$ system, and let $\{M_n\} \in D^{\leq 0}(\Fun(\N^\opp,\Mod_A))$ be another system. Then $\{K_n \otimes_A M_n\}$ is strict-essentially $0$, and
\[ R\lim_n (K_n \otimes_A M_n) \simeq 0.\]
\end{lemma}
\begin{proof}
The Milnor sequence takes the shape
\[ 1 \to \lim_n H^{-i}(K_n \otimes_A M_n) \to H^{-i}(R\lim_n (K_n \otimes_A M_n)) \to \lim^1_n H^{-i-1}(K_n \otimes_A M_n) \to 1.\]
Since both $\lim$ and $\lim^1$ vanish for a strict-essentially $0$ system of $A$-modules, it suffices to show that $\{H^{-i}(K_n \otimes_A M_n)\}$ is a strict-essentially $0$ system of $A$-modules for each $i \geq 0$. This follows immediately from the Kunneth spectral sequence which gives a (functorial in $n$) finite filtration on $H^{-i}(K_n \otimes_A M_n)$ whose graded pieces are subquotients of $\Tor^A_j(H^{-k}(K_n),M_n)$ for $j+k = i$ with $0 \leq j,k \leq i$.
\end{proof}

The following lemma plays a key role in relating derived and underived constructions later in this paper:

\begin{lemma}[Quillen]
\label{lem:sezeronoetheriantor}
Let $I \subset A$ be an ideal in a noetherian ring $A$, and let $M$ be a finitely generated $A$-module. Then the cone of the map $\{M \otimes_A A/I^n\} \to \{M/I^nM\}$ of objects in $D(\Fun(\N^\opp,\Mod_A))$ is strict-essentially $0$.
\end{lemma}
\begin{proof}
It suffices to check that $\{\Tor_i^A(M,A/I^n)\}$ is strict-essentially $0$ for $i > 0$.  By dimension shifting, it suffices to verify the claim for $i = 1$. Choose a presentation
\[ 1 \to R \to F \to M \to 1 \]
with $F$ a finite free $A$-module, and $R$ the kernel. Then 
\[ \Tor_1^A(M,A/I^n) \simeq \ker(R/I^nR \to F/I^nF) \simeq (I^n F \cap R)/I^n R,\]
so the claim follows from Example \ref{ex:artinrees}.
\end{proof}

\begin{corollary}
\label{cor:sezerocompletions}
Let $I \subset A$ be an ideal in a noetherian ring, and let $M$ be a finitely generated $A$-module. For any object $\{K_n\} \in D^{\leq 0}(\Fun(\N^\opp,\Mod_A))$, the natural map induces an equivalence
\[ \phi:R\lim_n(M \otimes_A A/I^n \otimes_A K_n) \simeq R\lim_n(M/I^nM \otimes_A K_n). \]
In particular, $\R\lim_n (M \otimes_A A/I^n) \simeq \widehat{M} := \lim_n M/I^nM$.
\end{corollary}
\begin{proof}
Consider the composite functor $F:D(\Fun(\N^\opp,\Mod_A)) \to D(\Fun(\N^\opp,\Mod_A)) \to D(\Mod_A)$ where the first one is $\{- \otimes_A K_n\}$, and the second one is $R\lim$. Then $F$ is an exact functor, and the map $\phi$ is simply $F$ applied to the natural map $\{M \otimes_A A/I^n\} \to \{M/I^n M\}$ in $D(\Fun(\N^\opp,\Mod_A))$. The claim now follows from Lemma \ref{lem:sezeronoetheriantor} and Lemma \ref{lem:sezerotensorproduct}.
\end{proof}

\begin{remark}
The object $\{M \otimes_A A/I^n\}$ referred to in Lemma \ref{lem:sezeronoetheriantor} and Corollary \ref{cor:sezerocompletions} is formally defined as $\{K \otimes_A A/I^n\} \in D(\Fun(\N^\opp,\Mod_A))$ for some $A$-flat resolution $K$ of $M$, and is independent of the choice of $K$.
\end{remark}

\begin{lemma}
\label{lem:adamscompeasy}
Let $f:A \to B$ be a ring map, and let $M$ be a $B$-module viewed as an $A$-module via $f$. Then the natural map $M \to M \otimes_A \Cech(f)$ is a homotopy-equivalence of cosimplicial $A$-modules. In particular, $M \simeq \Comp_A(M,f)$.
\end{lemma}
\begin{proof}
The $B$-action on $M$ defines the desired homotopy.
\end{proof}

\begin{proposition}[{\cite[Theorem 4.4]{CarlssonDercomp}}]
\label{prop:adamscompusual}
Let $I \subset A$ be an ideal in a noetherian ring $A$, and let $M$ be a finitely generated $A$-module. Then there is a natural isomorphism $\widehat{M} \simeq \Comp_A(M,I)$.
\end{proposition}
\begin{proof}
Let $F:D(\Mod_A) \to D(\Mod_A)$ be the exact functor $M \mapsto \Tot(M \otimes_A \Cech(f))$. There is an evident natural transformation $\eta:\id \to F$, i.e., a functorial map $\eta_M:M \to F(M)$. We first claim that $\eta_M$ is an equivalence when $M$ is an $A/I^n$-module for any integer $n \geq 1$. As both $\id$ and $F$ are exact functors, devissage reduces the claim to the case $n = 1$, where it follows from Lemma \ref{lem:adamscompeasy}. Hence, for any finitely generated $A$-module $M$, we have equivalences
\[ \eta_{M/I^nM}:M/I^nM \simeq \Tot(M/I^nM \otimes_A \Cech(f)).\]
Taking a limit over $n \in \N^\opp$ and commuting $R\lim$ with $\Tot$ then gives an equivalence
\[ \widehat{\eta}:\widehat{M} := R\lim_n M/I^nM \simeq \Tot(R\lim_n (M/I^nM \otimes_A \Cech(f))).\]
There is a natural map $\Tot(M \otimes_A \Cech(f)) \to \Tot(R\lim_n (M/I^nM \otimes_A \Cech(f)))$, so it suffices to check that $M \otimes_A \Cech(f)$ and $R\lim_n (M/I^nM \otimes_A \Cech(f))$ are equivalent in $D(\Fun(\Delta,\Mod_A))$ under the natural map 
\[\phi:M \otimes_A \Cech(f) \to R\lim_n (M/I^nM \otimes_A \Cech(f)).\]
The term in cosimplicial level $[m] \in \Delta$ in the source of $\phi$ is $M \otimes_A (A/I)^{\otimes (m+1)}$. The corresponding term on the target of $\phi$ is  $R\lim_n (M/I^n M \otimes_A (A/I)^{\otimes (m+1)})$.
Using Corollary \ref{cor:sezerocompletions} three times, we can rewrite this as
\begin{eqnarray*}
 R\lim_n (M/I^n M \otimes_A (A/I)^{\otimes (m+1)}) &=&  R\lim_n(M \otimes_A A/I^n \otimes_A (A/I)^{\otimes (m+1)})  \\
													 &=&  R\lim_n(M \otimes_A A/I^n \otimes_A A/I \otimes_A (A/I)^{\otimes m}) \\
													 &=&  R\lim_n(M \otimes_A (A/I)/I^n \otimes_A (A/I)^{\otimes m}) \\
													 &=&  R\lim_n(M \otimes_A (A/I)^{\otimes (m+1)}) \\
													 &=& M \otimes_A (A/I)^{\otimes (m+1)},
\end{eqnarray*}
which is exactly the source of $\phi$, proving that $\phi$ is an equivalence.
\end{proof}

\begin{remark}
If $A$ is any ring and $I = (f_1,\dots,f_r)$ is a finitely generated ideal, then, using the Postnikov filtration on a cosimplicial chain complex, one can show that $\Comp_A(A,I)$ coincides with the Greenlees-May \cite{GreenleesMay} derived $I$-adic completion (which is also studied by Lurie in \cite[\S 4]{LurieFormalExistence} for $E_\infty$-rings). 
\end{remark}

\section{Derived de Rham cohomology and the Adams completion}
\label{sec:mainthm}

Our goal in this section is to introduce the Hodge-completed version of Illusie's derived de Rham complex, and compare with the complex defined by Hartshorne using the Adams completion from \S \ref{sec:adamscomp}.
\subsection{Statement of the main theorem}

We first recall the definition of Illusie's derived de Rham cohomology, completed for the Hodge filtration.

\begin{construction}
\label{cons:ddr}
Let $A \to B$ be a map in $s\Alg$, and let $P \to B$ be a polynomial $A$-algebra resolution of $B$. Then we define the {\em Hodge-completed derived de Rham complex} $\widehat{\dR}_{B/A}$ of $A \to B$ as the completion of $|\Omega^*_{P/A}|$ for its Hodge filtration, i.e., $\widehat{\dR}_{B/A}$ is the limit over $\N^\opp$ of the diagram $\{\widehat{\dR}_{B/A}/\Fil^k_H\} \in \Fun(\N^\opp, \Ch(A))$ where $\widehat{\dR}_{B/A}/\Fil^k_H \in \Ch(A)$ is the totalisation of the simplicial cochain complex $[n] \mapsto \sigma^{\leq k} \Omega^*_{P/A}$ (where $\sigma^{\leq k}$ refers to the stupid truncation in cohomological degrees $\leq k$). This construction is independent of the choice of $P$ (up to homotopy), and has graded pieces computed by
\[ \gr^k_H(\widehat{\dR}_{B/A}) \simeq \wedge^k L_{B/A}[-k].\]
The above discussion also applies to any map of simplicial commutative rings in a topos. For a map $f:X \to Y$ of schemes, we write $\widehat{\dR}_{X/Y}$ for $\widehat{\dR}_{\calO_X/f^{-1} \calO_Y}$.
\end{construction}

\begin{remark}
A more explicit description of Construction \ref{cons:ddr} is: the complex $\widehat{\dR}_{B/A}$ is the completion of the second quadrant bicomplex $\Omega^*_{P/A}$ (where $\ast$ denotes the vertical variable) for the filtration defined by the rows. The homotopy-limit definition is more convenient if one wants to see the $E_\infty$-algebra structure and also work with diagrams of rings.
\end{remark}

\begin{remark}
	The associations $(A \to B) \mapsto \widehat{\dR}_{B/A}/\Fil^n_H$ and $(A \to B) \mapsto \widehat{\dR}_{B/A}$ from Construction \ref{cons:ddr} can be made strictly functorial by using the canonical free $A$-algebra resolution $P \to B$. In particular, given a diagram $F:I \to s\Alg_{A/}$ indexed by a small category $I$, there are induced diagrams $\widehat{\dR}_{F/A},\widehat{\dR}_{F/A}/\Fil^n_H:I \to \Ch(A)$ defined by $i \mapsto \widehat{\dR}_{F(i)/A}$ and $i \mapsto \widehat{\dR}_{F(i)/A}/\Fil^n_H$ (functorial in $F$), and an identification $\widehat{\dR}_{F/A} \simeq R\lim_n \widehat{\dR}_{F/A}/\Fil^n_H$ in $\Fun(I,\Ch(A))$; a morphism $F \to G$ of diagrams lets us define $\widehat{\dR}_{G/F}$ in the obvious way. Given such an $F$, taking limits over $I$ and commuting limits with limits gives us complexes $R\lim_I \widehat{\dR}_{F/A}/\Fil^n_H$ and $R\lim_I \widehat{\dR}_{F/A} \simeq R\lim_n R\lim_I \widehat{\dR}_{F/A}/\Fil^n_H$ in $\Ch(A)$; these objects will be used in the sequel for $I = \Delta$.
\end{remark}

We discuss a few examples of the derived de Rham complex:

\begin{example}
If $A \to B$ is a smooth morphism, then $L_{B/A} \simeq \Omega^1_{B/A}$, so $\widehat{\dR}_{B/A}$ is the usual de Rham complex.
\end{example}

\begin{example}
\label{ex:ddrlciquot}
Let $A$ be a $\Q$-algebra, and let $B = A/(f)$ where $f \in A$ is a regular element. Then the transitivity triangle degenerates to give an isomorphism $L_{B/A} \simeq (f)/(f^2)[1]$. In particular, $\widehat{\dR}_{B/A}/\Fil^2_H$ is an extension of $B$ by $L_{B/A}[-1] \simeq (f)/(f^2)$. One can check that the natural map $A \to \widehat{\dR}_{B/A}$ induces an equivalence $A/(f^2) \simeq \widehat{\dR}_{B/A}/\Fil^2_H$ (either via explicit simplicial resolutions, or by observing that the complex $\widehat{\dR}_{B/A}/\Fil^2_H := \Big(B \to L_{B/A}\Big)$ is the universal square zero extension of $B$ relative to $A$). In fact, one can show that the map $A \to \widehat{\dR}_{B/A}$ induces equivalences $A/(f^n) \simeq \widehat{\dR}_{B/A}/\Fil^n_H$, and hence a limiting equivalence $\widehat{A} \simeq \widehat{\dR}_{B/A}$. By Kunneth, one has $A/I^n \simeq \widehat{\dR}_{B/A}/\Fil^n_H$ for $B = A/I$ with $I \subset A$ any regular ideal.
\end{example}

\begin{example}
	Let $k$ be a field, and let $A$ be a finite type $k$-algebra which is not lci; for example, we may take $A = k[x,y]/(x^2,xy,y^2)$. Then $L_{A/k}$ is a left-unbounded complex by \cite{AvramovQuillenConj}, so $\widehat{\dR}_{A/k}/\Fil^2_H$ is also left-unbounded. In contrast, Theorem \ref{thm:ddralgdrcomp} shows that $\widehat{\dR}_{A/k} \simeq k$ when $k$ has characteristic $0$.
\end{example}

Next, we introduce Hartshorne's algebraic de Rham complex in the affine setting:

\begin{construction}
\label{cons:hartsalgdraff}
Let $f:A \to B$ be a finite type map of noetherian $\Q$-algebras, and fix a presentation $F \to B$ with $F$ a finite type polynomial $A$-algebra. Then we define the {\em algebraic de Rham complex} $\Omega^{H}_{B/A} \in D(\Mod_A)$ via
\[ \Omega^H_{B/A} := \Omega^*_{F/A} \otimes_F \widehat{F},\]
where $\widehat{F}$ denotes the completion of $F$ along the $I = \ker(F \to A)$; this complex is independent of choice of $F$. 
The construction  endows $\Omega^H_{B/A}$ with two filtrations: the filtration defined by the Hodge filtration on $\Omega^*_{F/A}$ is called the {\em formal} Hodge filtration (and depends on $F$), while the one obtained by tensoring the $I$-adic filtration on $\widehat{F}$ with the Hodge filtration on $\Omega^*_{F/A}$  is called the {\em infinitesimal} Hodge filtration (and is independent of $F$). The latter filtration is denoted by $\Fil^*_\inf$ (on both $\Omega^H_{B/A}$ and its cohomology), and is explicitly defined by
\[ \Omega^H_{B/A}/\Fil^p_\inf := \Big(F/I^p \to F/I^{p-1} \otimes_F \Omega^1_{F/A} \to F/I^{p-2} \otimes_F \Omega^2_{F/A} \to \dots \Big)\]
with the convention that $I^k = F$ for $k \leq 0$.
\end{construction}

\begin{remark}
\label{rmk:infderhamcomp}
Let $X = \Spec(B)$, $Y = \Spec(A)$, and $f:X \to Y$ be the map from Construction \ref{cons:hartsalgdraff}. The complex $\Omega^H_{B/A}$ is more conceptually the (derived) pushforward of the structure sheaf along $u:\Shv( (X/Y)_\inf) \to \Shv(X_\zar)$ defined by $(u^{-1} \calF)(U \hookrightarrow T) = \calF(U)$; the infinitesimal Hodge filtration on $\Omega^H_{B/A}$ corresponds to the filtration by powers of the defining ideal on the infinitesimal site (see \cite[Remark 3.7]{BhattdeJong} for an explanation). This description immediately shows that $\Omega^H_{B/A}$ is independent of $F$, and also extends naturally to the global setting discussed in \S \ref{ss:glob}.
 \end{remark}

\begin{example}
Let $A$ be a finite type $\C$-algebra. Then a theorem of Hartshorne shows that $\Omega^H_{A/\C}$ computes the Betti cohomology of $\Spec(A)^\an$. For example, if $A = \C[x,y]/(y^2 - x^3)$, then we may take $F = \C[x,y]$, so
\[ \Omega^H_{A/\C} \simeq \Big(\widehat{F} \to \widehat{F} dx \oplus \widehat{F} dy \to \widehat{F} (dx \wedge dy) \Big)\] 
where $\widehat{F}$ is the completion of $F$ along $(y^2 - x^3)$. Then $\Spec(A)^\an$ is contractible, so $R\Gamma(\Spec(A)^\an,\C) \simeq \C[0]$. It is a pleasant exercise to check that $\Omega^H_{A/\C}$ is also equivalent to $\C[0]$; see also \S \ref{ss:explicitex}.
\end{example}

The (affine version of our) main theorem is:

\begin{theorem}
\label{thm:ddralgdrcomp}
Let $A \to B$ be a finite type map of noetherian $\Q$-algebras. There is a filtered $A$-algebra 
\[ \widehat{\dR}_{B/A} \to \Omega^H_{B/A}\]
which is an equivalence of the underlying algebras.
\end{theorem}

We first prove a variant of Theorem \ref{thm:ddralgdrcomp} for certain maps of simplicial commutative rings using a convergence theorem of Quillen in \S \ref{ss:mainthmquillen}. In \S \ref{ss:mainthmgeneral}, we prove the theorem in general, and use it to give an elementary description of $\Omega^H_{B/A}$ promised in \S \ref{sec:intro}. This elementary description is reproven directly in \S \ref{ss:explicitex} in the case of graded rings.

\subsection{Quillen's convergence theorem: a review}
\label{ss:mainthmquillen}

We first recall Quillen's original theorem (see \cite[Theorem 8.8]{QuillenCRCNotes}), and then give a reformultion closer to derived de Rham theory.

\begin{proposition}[Quillen]
\label{prop:quillenconv}
Let $A \in s\Alg$, and let $I \subset A$ be a simplicial ideal with $\pi_0(I) = 0$. Assume that $I_n \subset A_n$ is a regular ideal. Then $A$ is $I$-adically complete, i.e., the natural map $A \simeq R\lim_n A/I^n$ is an equivalence in $s\Alg_{A/}$.
\end{proposition}
\begin{proof}
All ideals, modules and constructions that occur below are simplicial $A$-modules unless otherwise specified. By the convergence of the Postnikov filtration, it suffices to show that $I^{n+1}$ is $n$-connected, which we show by induction on $n$. Tensoring 
\[ 1 \to I \to A \to A/I \to 1\]
with $I^n$ and peeling of a suitable chunk gives an exact triangle
\[ \Big(\tau_{\geq 1} A/I \otimes_A I^n \Big) [-1] \to I \otimes_A I^{n} \to I^{n+1}\]
of simplicial $A$-modules. Note that $I$ is connected, and $I^{n}$ is $(n-1)$-connected by induction. The Kunneth spectral sequence shows that $I \otimes_A I^{n}$ is $n$-connected. Hence,  it suffices to show that $\Big(\tau_{\geq 1} A/I \otimes_A I^n \Big) [-1]$ is $n$-connected. This complex admits a complete decreasing separated exhaustive filtration\footnote{This is simply a way of expressing some qualitative features of one of the spectral sequences for a bisimplicial abelian group.} with graded pieces given by the simplicial $A/I$-complex $\Tor^A_i(A/I,I^n)[i-1]$ for $i \geq 1$ , so it suffices to show $n$-connectivity for each of these complexes. By dimension shifting, we have
\[ \Tor_{i+1}^A(A/I,A/I^n)[i-1] \simeq \Tor^A_{i}(A/I,I^n)[i-1],\]
so it suffices to show the $n$-connectivity of $\Tor_j^A(A/I,A/I^n)[j-2]$ for $j \geq 2$.  The regularity of the ideal $I$ gives a Koszul calculation of $\Tor^*_A(A/I,A/I^n)$ in each simplicial degree, and taking a homotopy colimit over $\Delta^\opp$ gives an exact sequence of simplicial $A/I$-modules
\[ 1 \to \wedge^{j+n}(I/I^2) \to \dots \to \wedge^{j+1}(I/I^2) \otimes_{A/I} \Sym^{n-1}(I/I^2) \to \wedge^j(I/I^2) \otimes_{A/I} \Sym^n(I/I^2) \to \Tor_j^A(A/I,A/I^n) \to 1.\]
Since $I$ is connected with $I/I^2$ termwise flat over $A/I$, the simplicial $A/I$-module $\wedge^a(I/I^2) \otimes_{A/I} \Sym^b(I/I^2)$ is $(a+b-2)$-connected for any $a,b \geq 0$. Hence, the preceding exact sequence together with the preservation of $k$-connectivity under homotopy-colimits shows that $\Tor_j^A(A/I,I^n)$ is $(n+j-2)$-connected. Hence, $\Tor_j^A(A/I,A/I^n)[j-2]$ is $n$-connected for $j \geq 2$, as wanted.
\end{proof}

\begin{remark}
If $I \subset A$ is an ideal in a simplicial commutative ring $A$, then the inverse system $\{\pi_0(A/I^n)\}$ has surjective transition maps, and hence is acyclic for $\lim$. As the functor of $\N^\opp$-indexed limits  has cohomological dimension $1$, it follows that the limit of $\{A/I^n\}$ in $s\Alg_{A/}$ and $D(\Mod_A)$ is the same. In particular, the limit appearing in Proposition \ref{prop:quillenconv} can be also be computed in $D(\Mod_A)$.
\end{remark}

\begin{remark}
\label{rmk:quillenconvhodge}
Let $I \subset A$ be as in Proposition \ref{prop:quillenconv}. Then the $I$-adic filtration $\Fil^*_H$ on $A$ is a decreasing complete separated exhaustive filtration by Proposition \ref{prop:quillenconv}. We can also identify the graded pieces as follows: by regularity of $I$, there is an equivalence  equivalence
\[ \gr^k_H(A) \simeq I^k/I^{k+1} \simeq \Sym^k(I/I^2)\] 
of simplicial $A/I$-modules. If $\Q \subset A$, then we can identify $\Gamma^k$ with $\Sym^k$ to write
\[ \gr^k_H(A) \simeq \Gamma^k(I/I^2).\] 
Using regularity of $I$ (and the transitivity triangle), there is an identification $I/I^2[1] \simeq L_{(A/I)/A}$, so we obtain
\[ \gr^k_H(A) \simeq \wedge^k L_{(A/I)/A}[-k].\]
Since the cotangent complex and its wedge powers are intrinsic to the map $A \to A/I$ in $s\Alg_{\Q/}$, one can then show that the filtration $\Fil^*_H$ on $A$ is intrinsic to the map $A \to A/I$ in $\Ho(s\Alg_{\Q/})$.  Going one step further, under the assumption $\Q \subset A$, one can identify $A_n/I_n^k \simeq \widehat{\dR}_{(A_n/I_n)/A_n}/\Fil^k_H$ for each $n$ using the regularity of $I_n \subset A_n$ (see Example \ref{ex:ddrlciquot}). Taking a colimit over $\Delta^\opp$ gives an identification $\{A/\Fil^k_H\} \simeq \{\widehat{\dR}_{(A/I)/A}/\Fil^k_H\}$ of $\N^\opp$-indexed diagrams. Hence, there is a limiting equivalence $A \simeq \widehat{\dR}_{(A/I)/A}$ for all simplicial $\Q$-algebras $A$ equipped with a termwise regular ideal $I \subset A$ satisfying $\pi_0(I) = 0$.
\end{remark}

We can now give the promised reformulation of Quillen's theorem; in essence, the last statement of Remark \ref{rmk:quillenconvhodge} is true without the regularity assumption. This result can be viewed as a version of Theorem \ref{thm:ddralgdrcomp} for maps between simplicial $\Q$-algebras with a connected kernel.

\begin{corollary}
\label{cor:quillenconvnoreg}
Let $A \in s\Alg_{\Q/}$, and let $I \subset A$ be an ideal with $\pi_0(I) = 0$. Then $A$ admits a functorial complete separated $\N^\opp$-indexed filtration $\Fil^k_H$ whose associated $\N^\opp$-indexed system of quotients is identified as
\[ \{A/\Fil^k_H\} \simeq \{\widehat{\dR}_{(A/I)/A}/\Fil^k_H\}.\]
In particular, one has $A \simeq \widehat{\dR}_{(A/I)/A}$.
\end{corollary}
\begin{proof}
Using \cite[\S II.4, Theorem 4]{QuillenHA}, one can define a simplicial model category structure on the simplicial category of pairs $(A,I)$ comprising $A \in s\Alg$ together with simplicial ideals $I \subset A$ as follows: a map $(A,I) \to (B,J)$ is a (trivial) fibration if and only if $A \to B$ and $I \to J$ are so as maps of simplicial sets. The cofibrant objects are pairs $(F,J)$ with each $F_n$ a polynomial algebra on a set $X_n$ of generators, and each $J_n \subset F_n$ an ideal defined by a subset $Y_n \subset X_n$ of the polynomial generators such that both $X_n$ and $Y_n$ are preserved by the degeneracies (and a similar relative version). In particular, for each cofibrant object $(F,J)$, the ideal $J$ is (termwise) regular. The desired claim now follows from the arguments of Proposition \ref{prop:quillenconv} and Remark \ref{rmk:quillenconvhodge} by passage to cofibrant replacements. 
\end{proof}

\begin{remark}
The proof of Corollary \ref{cor:quillenconvnoreg} shows that the Hodge filtration $\Fil^k_H$ on $A$ is {\em finer} than the $I$-adic filtration, i.e., there is a map $\{A/\Fil^k_H\} \to \{A/I^k\}$ of $\N^\opp$-indexed systems inducing an equivalence on $R\lim$, but not necessarily on $\gr^k_H$.
\end{remark}

\subsection{The main theorem}
\label{ss:mainthmgeneral}

The goal of this section is to prove Theorem \ref{thm:ddralgdrcomp} and record some corollaries. We begin with the case of surjective maps of rings, which is the heart of Theorem \ref{thm:ddralgdrcomp}. In fact, we prove a more general result for surjective maps in $s\Alg$, i.e., maps in $s\Alg$ which are surjective on $\pi_0$.

\begin{proposition}
\label{prop:completionfiltration}
Let $f:A \to B$ be a surjective map in $s\Alg_{\Q/}$.  Then $\Comp_A(A,f) \in \Ch(A)$ admits a canonical $\N^\opp$-indexed separated complete filtration $\Fil^\bullet_H$ whose associated $\N^\opp$-indexed system of quotients is identified as
\[ \{\Comp_A(A,f)/\Fil^k_H\} \simeq \{\widehat{\dR}_{B/A}/\Fil^k_H \}.\]
In particular, one has $\Comp_A(A,f) \simeq \widehat{\dR}_{B/A}$.
\end{proposition}

\begin{proof}
We may assume $B$ is $A$-cofibrant, i.e., $B$ is a simplicial polynomial $A$-algebra. By definition, one has $B_m := \Cech(A \to B)([m]) \simeq B^{\otimes_A [m]}$ with the augmentation $\Cech(A \to B) \to B$ corresponding to the multiplication map. Since $A \to B$ is surjective, one has $\pi_0(B_m) \simeq \pi_0(B)$. Hence, by Corollary \ref{cor:quillenconvnoreg}, there is a functorial complete separated $\N^\opp$-indexed filtration $\Fil^*_H$ on $B_m$ with an identification of $\N^\opp$-indexed systems $\{B_m/\Fil^k_H\} \simeq \{ \widehat{\dR}_{B/B_m}/\Fil^k_H \}$. By functoriality, these filtrations fit together to define a functorial complete separated $\N^\opp$-indexed filtration $\Fil^*_H$ on $\Cech(A \to B)$ together with identifications of $\N^\opp$-indexed systems of cosimplicial $\Cech(A \to B)$-complexes
\[  \{\Cech(A \to B)/\Fil^k_H\}  \simeq \{\widehat{\dR}_{B/\Cech(A \to B)}/\Fil^k_H\}. \]
As limits commutes with limits, there is an induced functorial complete separated $\N^\opp$-indexed filtration $\Fil^*_H$ on $\Comp_A(A,f)$ with an identification $\N^\opp$-indexed systems
\[ \{\Comp_A(A,f)/\Fil^k_H\} \simeq \{ \Tot(\widehat{\dR}_{B/\Cech(A \to B)}/\Fil^k_H) \}.\]
On the other hand, functoriality induces a natural map of $\N^\opp$-indexed systems 
\[ \epsilon:\{\widehat{\dR}_{B/A}/\Fil^k_H\} \to \{ \Tot(\widehat{\dR}_{B/\Cech(A \to B)}/\Fil^k_H) \},\] 
so it suffices to check that $\epsilon$ is an equivalence. Passing to graded pieces, it suffices to check that the natural map
\[ \gr^k(\epsilon):\wedge^k L_{B/A}[-k] \to \Tot(\wedge^k L_{B/\Cech(A \to B)}[-k]),\]
which is proven in Corollary \ref{cor:cotcomplexcosimpalg}.
\end{proof}

Proposition \ref{prop:completionfiltration} shows the ``topological invariance'' of the $\widehat{\dR}_{-/k}$ functor:

\begin{example}
	Fix a $\Q$-algebra $k$, and let $f:k \to A$ be a {\em connected} simplicial $k$-algebra, i.e., $k \simeq \pi_0(A)$. For example, we can take $A = \Sym_k(K)$ for any complex $K \in D^{< 0}(\Mod_k)$. Then $f$ is surjective, so Proposition \ref{prop:completionfiltration} implies $\Comp_A(k,f) \simeq \widehat{\dR}_{A/k}$. On the other hand, since $f$ has a section (via the truncation map $A \to \pi_0(A) \simeq k$),  Example \ref{ex:cechnervesection} shows $k \simeq \Comp_A(k,f)$. In particular, $k \simeq \widehat{\dR}_{A/k}$. 
\end{example}

\begin{example}
	Let $f:A \to B$ be a surjective map of noetherian $\Q$-algebras with nilpotent kernel. Then Proposition \ref{prop:completionfiltration} shows that $\Comp_A(A,f) \simeq \widehat{\dR}_{B/A}$. On the other hand, since $A$ is noetherian, Proposition \ref{prop:adamscompusual} shows $\Comp_A(A,f) \simeq A$, so $A \simeq \widehat{\dR}_{B/A}$.
\end{example}

We can now prove Theorem \ref{thm:ddralgdrcomp} by bootstrapping from Proposition \ref{prop:completionfiltration}.

\begin{proof}[Proof of Theorem \ref{thm:ddralgdrcomp}]
Choose a finite type polynomial $A$-algebra $F$ and an $A$-algebra surjection $F \to B$. Then $\Omega^H_{B/A} \simeq \Omega^*_{F/A} \otimes_F \widehat{F}$. It is well-known (see \cite{BhattdeJong}, for example) that
\[ \Omega^H_{B/A} \simeq \Tot(\cosimp{\widehat{F}}{\widehat{F \otimes_A F}}{\widehat{F \otimes_A F \otimes_A F}})\]
where the completion takes place along the composition $F^{\otimes_A [n]} \to F \to B$. By Proposition \ref{prop:adamscompusual}, we can write
 \[ \Omega^H_{B/A} \simeq \Tot(\Comp_A(\Cech(A \to F) \to B)).\]
As the maps $F^{\otimes_A [n]} \to B$ are surjective, Proposition \ref{prop:completionfiltration} gives
\[ \Omega^H_{B/A} \simeq \Tot(\widehat{\dR}_{B/\Cech(A \to F)}).\]
This formula equips $\Omega^H_{B/A}$ with a complete separated $\N^\opp$-indexed filtration $\Fil^*_{H'}$ whose associated $\N^\opp$-indexed system of quotients is given by
\[ \{ \Omega^H_{B/A}/\Fil^k_{H'} \} \simeq \{\Tot(\widehat{\dR}_{B/\Cech(A \to F)}/\Fil^k_H)\}.\]
Repeating the argument at the end of the proof of Proposition \ref{prop:completionfiltration} then finishes the proof.
\end{proof}

Theorem \ref{thm:ddralgdrcomp} yields the following explicit representative for $\Omega^H_{B/A}$ which seems to be new:

\begin{corollary}
\label{cor:elemderham}
Let $f:A \to B$ be a finitely presented flat map of noetherian $\Q$-schemes. Then 
\[ \Omega^H_{B/A} \simeq \Tot(\widehat{\Cech}(A \to B)) := \Tot\Big(\cosimp{B}{\widehat{B \otimes_A B}}{\widehat{B \otimes_A B \otimes_A B}} \Big).\]
where the completions take place along the multiplication maps $B^{\otimes n} \to B$.
\end{corollary}
\begin{proof}
Let $K$ denote the complex (or $E_\infty$-ring) appearing on the right above. Using propositions \ref{prop:adamscompusual} and \ref{prop:completionfiltration} as in the proof of Theorem \ref{thm:ddralgdrcomp} as well as the flatness of $f$, the complex $K$ can be rewritten as
\[ K \simeq \Tot(\widehat{\dR}_{B/\Cech(A \to B)}).\]
In particular, there exists a functorial complete separated $\N^\opp$-indexed filtration $\Fil^*_H$ on $K$ with an identification of $\N^\opp$-indexed systems
\[ \{K/\Fil^k_H \} \simeq \{ \Tot(\widehat{\dR}_{B/\Cech(A \to B)}/\Fil^k_H)\}.\]
The argument at the end of the proof of Proposition \ref{prop:completionfiltration} then identifies the right hand side with $\{\widehat{\dR}_{B/A}/\Fil^k_H\}$. The claim then follows by taking limits and using Theorem \ref{thm:ddralgdrcomp}.
\end{proof}

\begin{remark}
The ring $\widehat{\Cech}(A \to B)$ occurring in Corollary \ref{cor:elemderham} is the completion of $\Cech(A \to B)$ along the augmentation map $\Cech(A \to B) \to B$. In particular, $\widehat{\Cech}(A \to B)$ admits a natural complete separated $\N^\opp$-indexed filtration (by powers of the augmentation ideal). We do not know what the induced filtration on $\Omega^H_{B/A}$ means.
\end{remark}

\begin{remark}
The equivalence $\Omega^H_{B/A} \simeq \Tot(\widehat{\Cech}(A \to B))$ from Corollary \ref{cor:elemderham} is also valid for non-flat maps if one interprets $\widehat{\Cech}(A \to B)$ to be the Adams completion of the (termwise) surjective map $\Cech(A \to B) \to B$. In fact, when stated as such,  one need not impose any finiteness conditions on the map $A \to B$. For example, if $A \to B$ is itself surjective, then 
\[\Comp_A(\Cech(A \to B) \to B) \simeq \Cech(A \to B)\]
as $\Comp_A(B^{\otimes [m]} \to B) \simeq B^{\otimes [m]}$ by \cite[Theorem 6.1]{CarlssonDercomp}; the claim now follows Proposition \ref{prop:completionfiltration}
\end{remark}

\begin{remark}
	The completion is absolutely essential in Corollary \ref{cor:elemderham}. Indeed, when $A \to B$ is faithfully flat, $\Tot(\Cech(A \to B))$ is simply $A$ by Example \ref{ex:cechnervefppf}, which is certainly not the value of relative de Rham cohomology in general. For an explicit example to ``see'' the effect of completion, consider the case where $A \to B$ is a finite \'etale cover. Then the cosimplicial ring $\Cech(A \to B)$ totalises to $A$ as before. On the other hand, the cosimplicial ring $\widehat{\Cech}(A \to B)$ is the constant cosimplicial ring with value $B$ as the completion of $B^{\otimes_A [m]} \to B$ along the (finite \'etale) multiplication map is simply $B$. Hence, $\Tot(\widehat{\Cech}(A \to B)) \simeq B$; this is the desired answer since $\Omega^*_{B/A} = B[0]$.
\end{remark}

\begin{example}
Let $R$ be a finite type $\C$-algebra. Then Corollary \ref{cor:elemderham} gives the following elementary description of the Betti cohomology of the associated complex analytic space:
\[ \R\Gamma(\Spec(R)^{\an},\C) \simeq \Tot\Big(\cosimp{R}{\widehat{R \otimes_\C R}}{\widehat{R \otimes_\C R \otimes_\C R}}\Big).\]
When $R$ is smooth, the above description is due to Grothendieck (through the infinitesimal site). At the opposite extreme, if $R$ is an Artinian local $\C$-algebra, then $R^{\otimes n}$ is complete for the filtration defined by any ideal, so $\C \simeq \Tot(\Cech(\C \to R)) \simeq \Tot(\widehat{\Cech}(\C \to R))$, as expected (since $\Spec(R)$ is homeorphic to a point); some further examples are discussed in \S \ref{ss:explicitex}.  Another corollary of this description is:  the complex on the right above is acyclic outside cohomological degrees in $[0,\dim(R)]$ by Artin's theorem on the cohomological dimension of affine schemes.
\end{example}

\begin{remark}
\label{rmk:poscharamitsur}
We do not know a good description for the complex $\Tot(\widehat{\Cech}(A \to B))$ in characteristic $p$ in general.  However, a variant operation does have meaning: the category of pairs from the proof of Proposition \ref{cor:quillenconvnoreg} supports a {\em derived} pd-completion functor obtained by deriving Berthelot's (classical) completed pd-envelope. One can then show that the derived pd-completion of $\Cech(A \to B)$ along its augmentation computes $\widehat{\dR}_{B/A}$ in any characteristic.
\end{remark}

\subsection{Globalisation}
\label{ss:glob}

To globalise Theorem \ref{thm:ddralgdrcomp}, we need the global version of Construction \ref{cons:hartsalgdraff}:

\begin{construction}
\label{cons:hartsalgdrgeom}
Let $f:X \to Y$ be a finite type morphism of noetherian $\Q$-schemes. Assume there exists an closed immersion $i:X \hookrightarrow Z$ with $Z$ smooth over $Y$. Set
\[ \Omega^H_{X/Y} = \widehat{i^{-1} \calO_Z} \otimes_{i^{-1} \calO_Z} i^{-1} (\Omega^*_{Z/Y}) \in D(\Mod_{f^{-1}\calO_Y}). \]
Here $\widehat{i^{-1} \calO_Z} = \lim_n i^{-1}(\calO_Z/I_X^n)$ is the completion of $i^{-1} \calO_Z$ along the ideal $i^{-1} (I_X)$ defining $\calO_X$; this construction is independent of choice of $Z$. We obtain the {\em formal} and {\em infinitesimal} Hodge filtrations as in Construction \ref{cons:hartsalgdraff}; denote the latter (on $\Omega^H_{X/Y}$ and its cohomology) by $\Fil^*_\inf$
\end{construction}

\begin{remark}
\label{rmk:localisationhartshorne}
In the notation of Construction \ref{cons:hartsalgdrgeom}, as cohomology commutes with (derived) limits, one has
\[ R\Gamma(X,\widehat{i^{-1} \calO_Z}) \simeq R\lim_n R\Gamma(X, i^{-1}(\calO_Z/I_X^n)) \simeq R\lim_nR\Gamma(Z,\calO_Z/I_X^n).\] 
In particular, if $Z$ (and hence $X$) is affine, then $R\Gamma(X,\widehat{i^{-1} \calO_Z})$ is simply the formal completion of the ring $\calO(Z)$ along the kernel of $\calO(Z) \to \calO(X)$. Similar reasoning shows that $R\Gamma(X,\Omega^H_{X/Y}) \simeq \Omega^H_{\calO(X)/\calO(Y)}$ if $X$, $Y$, and $Z$ are affine, i.e., Construction \ref{cons:hartsalgdrgeom} is compatible with Construction \ref{cons:hartsalgdraff}.
\end{remark}

Remark \ref{rmk:localisationhartshorne} explains that the Hartshorne's algebraic de Rham complex for a map of schemes recovers the affine one discussed in Construction \ref{cons:hartsalgdraff} when taking sections over affines. Likewise, the Hodge-completed derived de Rham complex also localises. Since the constructions of \S \ref{ss:mainthmgeneral} were functorial in the rings, they make sense for diagrams of rings, such as a cosimplicial one obtained from an affine open cover. Iterating this observation twice gives:

\begin{corollary}
\label{cor:ddrgeomcomp}
Let $f:X \to Y$ be a finite type map of noetherian $\Q$-schemes. Assume that $X$ can be realised as a closed subscheme of a smooth $Y$-scheme. Then there is a natural filtered $f^{-1}\calO_Y$-algebra map
\[  \widehat{\dR}_{X/Y} \to \Omega^H_{X/Y} \]
that is an equivalence of the underlying algebras. 
\end{corollary}

\begin{remark}
As explained in Remark \ref{rmk:infderhamcomp}, the complex $\Omega^H_{X/Y}$ from Corollary \ref{cor:ddrgeomcomp} may be identified with $\R u_* (\calO_{X/Y,\inf})$ where $u:\Shv((X/Y)_\inf) \to \Shv(X_\zar)$ is the structure map for infinitesimal site of $X$ relative to $Y$. Corollary \ref{cor:ddrgeomcomp} can be then reformulated to give an equivalence $\R u_* (\calO_{X/Y,\inf}) \simeq \widehat{\dR}_{X/Y}$. This statement makes no reference to the ambient smooth $Y$-scheme $Z$, and is indeed true without the assumption that $X$ admits such an embedding.
\end{remark}

\begin{remark}
The map in Corollary \ref{cor:ddrgeomcomp} is typically {\em not} a filtered equivalence. For example, if $A \to B = A/I$ is a surjective map of noetherian $\Q$-algebras, then the map $\widehat{\dR}_{B/A} \to \Omega^H_{B/A} \simeq \widehat{A}$ induces on $\gr^1$ the derivation $L_{B/A}[-1] \to I/I^2$ that classifies the square-zero extension $A/I^2 \to B$ of $B$ by $I/I^2$, and is an isomorphism only if $I$ is a regular ideal; the maps on higher graded pieces are computed using multiplicativity.
\end{remark}

\begin{corollary}
\label{cor:elemdrglob}
Let $f:X \to Y$ be a finite type flat map of noetherian $\Q$-schemes. Let $f^{(n)}:\widehat{X(n)} \to Y$ denote the formal completion of the $(n+1)$-fold fibre-product of $f$ along the diagonal. Then there is an equivalence
\[ \R f_* (\Omega^H_{X/Y}) \simeq \Tot\Big(\cosimp{\R f_* (\calO_X)}{\R f^{(1)}_* (\calO_{\widehat{X(1)}})}{\R f^{(2)}_* (\calO_{\widehat{X(2)}})} \Big). \]
\end{corollary}

\begin{remark}
In Corollary \ref{cor:elemdrglob}, if the map $f$ is smooth, then the result is essentially immediate once we know that (as explained in Remark \ref{rmk:infderhamcomp}) the complex $\R f_*(\Omega^H_{X/Y})$ can be identified with $\R f^\inf_* (\calO_{X/Y,\inf})$ where $f^\inf:(X/Y)_\inf \to Y$ is the structure map for infinitesimal site of $X$ relative to $Y$. Indeed, by the smoothness of $f$, the trivial thickening $X \hookrightarrow X$ defines an object of $(X/Y)_\inf$ dominating the final object of $\Shv( (X/Y)_\inf)$, and the claim follows by standard Cech theory. The surprising assertion in Corollary \ref{cor:elemdrglob} is that the same recipe also works for any flat map $f$, even though the object $X \hookrightarrow X$ certainly does not cover the final object of $\Shv( (X/Y)_\inf)$ in general.
\end{remark}

%\begin{remark}
%\label{rmk:derinfsite}
%We suspect Corollary \ref{cor:ddrgeomcomp} can be reformulated topologically. More precisely, there should exist a {\em derived} infinitesimal $\infty$-site $(X/Y)_{d\inf}$ parametrising derived thickenings of Zariski open subsets of $X$ relative to $Y$, and a map $\pi:\Shv( (X/Y)_\inf) \to \Shv( (X/Y)_{d\inf})$ with $\pi_\ast$ induced by the evident inclusion $(X/Y)_\inf \subset (X/Y)_{d\inf}$. Corollary \ref{cor:ddrgeomcomp} would then follow if we showed that $(X/Y)_\inf \subset (X/Y)_{d\inf}$ is cofinal, i.e., if every derived thickening can be dominated by a classical one. This should follow from Quillen's Nilpotent extension theorem \cite[Theorem 9.4]{QuillenCRCNotes}, but we do not develop this line of reasoning here.
%\end{remark}

\subsection{An explicit example}
\label{ss:explicitex}

Let $f:A \to B$ be a finite type flat map of noetherian $\Q$-algebras. Let $\widehat{\Cech}(A \to B)$ be the completion of $\Cech(A \to B)$ along the augmentation $\Cech(A \to B) \to B$. Corollary \ref{cor:elemderham} shows that $\Tot(\widehat{\Cech}(A \to B))$ computes the algebraic de Rham cohomology of $f$.  Here we reprove this result ``by hand'' in a special case: $A = \C$ and $B$ is a homogeneous singularity. The main idea is to prove homotopy-invariance for the functor $R \mapsto \Tot(\widehat{\Cech}(\C \to R))$, and then exploit the contracting $\A^1$-action on $\Spec(B)$.

We start with a general criterion for acyclicity of certain complexes.

\begin{lemma}
\label{lem:acycmod}
Let $A$ be a ring complete for the topology defined by some ideal $I \subset A$.  Let $K \in \Ch^{\geq 0}(A)$ be a cochain complex of flat $A$-modules with $I$-adically complete terms. If $K/IK$ is acyclic, so is $K$.
\end{lemma}
\begin{proof}
The flatness hypothesis and devissage give $K/I^nK \simeq 0$ for all $n$, so it suffices to prove $K \simeq R\lim_n K/I^nK$. As $K$ is a cochain complex, we can write $K \simeq R\lim_m K/\sigma^{\geq m} K$ where $\sigma^{\geq m}$ denotes the stupid truncation in cohomological degrees $\geq m$. Commuting limits, we reduce to the case where $K$ has only finitely many non-zero terms, where the claim is clear as $K$ has $I$-adically complete terms.
\end{proof}

\begin{remark}
	In Lemma \ref{lem:acycmod}, we do not know if  $K/IK$ agrees with the derived tensor product $K \otimes_A A/I$.
\end{remark}

Fix a ring $A$ with an ideal $I$. For any map $f:A \to B$, we write $\widehat{\Cech}_I(A \to B)$ for the (termwise) completion of $\Cech(A \to B)$ along the composite $\Cech(A \to B) \to B \to B/IB$. We have the following:

\begin{lemma}[Poincare Lemma]
\label{lem:poincare}
Let $A$ be noetherian $\Q$-algebra, and let $\widehat{A}$ be its completion for the topology defined by an ideal $I \subset A$. Then $\Tot(\widehat{\Cech}_I(A \to A[t])) \simeq \widehat{A}$.
\end{lemma}
\begin{proof}
Since reduction modulo $I$ commutes with completion, we may assume that $A = \widehat{A}$. Let $K$ denote the cochain complex underlying the augmented cosimplcial $A$-algebra $A \to \widehat{\Cech}_I(A \to A[t])$. Then $K$ has $A$-flat and $I$-adically complete terms. Moreover, since completion is exact (on finitely generated modules over noetherian rings), one also has an identification in $c\Alg_{A/}$
\[ \widehat{\Cech}_I(A \to A[t]) \otimes_A A/I \simeq \widehat{\Cech}(A/I \to A/I[t]).\]
The formal Poincare lemma and the characteristic $0$ hypothesis then show that the totalisation of the right hand side is simply $A/I$ via the natural map. It follows that $K/IK$ is acyclic. By Lemma \ref{lem:acycmod}, $K$ is acyclic too, as desired.
\end{proof}

Using the Poincare lemma, we can show (for a fixed characteristic $0$ field $\C$):

\begin{proposition}[Homotopy invariance]
\label{prop:homotopyinvariance}
Let $R$ be finite type $\C$-algebra. Then $R \to R[t]$ induces an equivalence
\[ \Tot(\widehat{\Cech}(\C \to R)) \simeq \Tot(\widehat{\Cech}(\C \to R[t])). \]
\end{proposition}
\begin{proof}
Let $A \in \Fun(\Delta \times \Delta,\Alg_{\C/})$ be the bicosimplicial ring obtained by taking tensor products of $\Cech(\C \to R)$ (viewed horizontally) with $\Cech(\C \to \C[t])$ (viewed vertically), and let $B$  denote the (termwise) completion of $A$ for its augmentation $A \to R \otimes_\C \C[t]$, i.e., we have
\[ B([m],-) = \widehat{\Cech}_I(R^{\otimes [m]} \to R^{\otimes [m]}[t]) \quad \mathrm{and} \quad B(-,[n]) = \widehat{\Cech}_J(\C[t]^{\otimes [n]} \to R[t]^{\otimes [n]}) \]
where $I \subset R^{\otimes [m]}$ is the kernel of the augmentation $R^{\otimes [m]} \to R$, and $J \subset \C[t]^{\otimes [n]}$ is the kernel of the augmentation $\C[t]^{\otimes [n]} \to \C[t]$.  It is easy to see that the diagonal cosimplicial ring $B|_{\Delta \subset \Delta \times \Delta}$ coincides with $\widehat{\Cech}(\C \to R[t])$. As the diagonal of a bicosimplicial abelian group computes its homotopy-limit, we get
\[ \Tot(\widehat{\Cech}(\C \to R[t])) \simeq R\lim_{\Delta \times \Delta} B.\]
On the other hand, computing the vertical totalisations first and using Lemma \ref{lem:poincare} shows that
\[ \widehat{\Cech}(\C \to R) \simeq \R\lim_\Delta B \in c\Alg_{\C/}\]
where $R\lim_\Delta$ denotes the homotopy-limit functor along the first projection $\Delta \times \Delta \to \Delta$. The claim now follows by taking limits over the residual $\Delta$ and comparing with the previous formula.
\end{proof}

Here is the promised example:

\begin{example}
Let $R$ be an $\N$-graded finite type $\C$-algebra with $R_0 = \C$. Then $\Spec(R)^{\an}$ is contractible, so $\R\Gamma(\Spec(R)^{\an},\C) \simeq \C$. We will show that $\Tot(\widehat{\Cech}(\C \to R)) \simeq \C$. The $\N$-grading defines an action of the multiplicative monoid scheme $\A^1$ on $\Spec(R)$ by dilations; explicitly, the grading gives a coaction map $R \to R[t]$ via
\[ \sum_{i \geq 0} a_i \mapsto \sum_{i \geq 0} a_i t^i \]
where $a_i \in R_i$. In particular, there is a family $\{f_t:\Spec(R) \to \Spec(R)\}$ of endomorphisms of $\Spec(R)$ parametrised by $t \in \C = \A^1(\C)$ such that $f_1 = \id_{\Spec(R)}$ while $f_0$ comes from the contraction $R \twoheadrightarrow R_0 = \C \hookrightarrow R$. By Proposition \ref{prop:homotopyinvariance}, pulling back along any $f_t$ induces the same map on cohomology. The claim now follows by comparing behaviours at $t = 0$ and $t = 1$.
\end{example}

\section{The derived Hodge filtration}
\label{sec:derhodgefilt}

For any algebraic variety $X$ over a characteristic $0$ field $k$, Corollary \ref{cor:ddrgeomcomp} gives rise to a spectral sequence
\[ E_1^{p,q}:H^q(X,\wedge^p L_{X/k}) \Rightarrow H^{p+q}(X,\Omega^H_{X/k})  \]
which induces the derived Hodge filtration $\Fil^*_H$ on the target.  Our goal in this section is to analyse this filtration and its relation to other better studied filtrations on the same group, such as the Hodge-Deligne filtration or the infinitesimal filtration.  First, we recall the Deligne-Du Bois complex \cite{DuBois} and the Hodge-Deligne filtration \cite{DeligneHodgeII}:

\begin{construction}
Let $X$ be a finite type $k$-scheme. The {\em Deligne-Du Bois} complex $\underline{\Omega}^*_{X/k}$ is defined as follows: fix a proper hypercover $f:X_\bullet \to X$ with each $X_n$ smooth (possible by Hironaka or de Jong), and set
\[ \underline{\Omega}^*_{X/k} := \Tot \big(\R f_* \Omega^*_{X_\bullet/k}\big) \quad \mathrm{in} \quad D^{\geq 0}(\Mod(X,k)).\]
The Hodge-Deligne filtration $\Fil^*_{HD}$ on $\underline{\Omega}^*_{X/k}$ refers to the object of the filtered derived category  $D^{\geq 0}(\Fun(\N^\opp,\Mod(X,k)))$ lifting $\underline{\Omega}^*_{X/k}$ defined by totalising the Hodge filtration on each $\Omega^*_{X_n/k}$. The graded pieces are computed by
\[ \gr^n_{HD} \underline{\Omega}^*_{X/k} \simeq \Tot \big( \R f_* \Omega^n_{X_\bullet/k} [-n]\big). \]
This filtered object is independent of the choice of $X_\bullet$, and may be characterised intrinsically as the $h$-localisation of the (filtered) presheaf of de Rham complexes on the site of all $X$-schemes \cite{LeedRhlocal}. If $k = \C$ and $X$ is proper, then the induced filtration on $H^*(X^\an,\C)$ is simply the Hodge filtration from mixed Hodge theory \cite{DeligneHodgeII}; we call it the Hodge-Deligne filtration here to distinguish it from the other filtrations.
\end{construction}

We now record the relation between certain natural filtrations:

\begin{proposition}
\label{prop:derhodgefilt}
Let $X$ be a finite type $k$-scheme. There are natural maps
\[ \widehat{\dR}_{X/k} \stackrel{a}{\to} \Omega^H_{X/k} \stackrel{b}{\to} \Omega^*_{X/k} \stackrel{c}{\to} \underline{\Omega}^*_{X/k} \]
of filtered complexes such that $a$, $c \circ b$, and $c \circ b \circ a$ induce an equivalence on the underlying complexes. In particular, the algebraic de Rham cohomology of $X$ is a summand of the cohomology of $\Omega^*_{X/k}$. 
\end{proposition}
\begin{proof}
	The filtered map $a$ is from Theorem \ref{thm:ddralgdrcomp} and an equivalence of the underlying complexes. The map $b$ comes from the functoriality of the Kahler differentials and can be easily checked to map the infinitesimal Hodge filtration on $\Omega^H_{X/k}$ (see Construction \ref{cons:hartsalgdraff}) to the stupid filtration on $\Omega^*_{X/k}$. The map $c$ also comes from the functoriality of the Kahler differentials and is obviously filtered. The rest of the assertions are immediate as $\Omega^H_{X/k}$ (and hence $\widehat{\dR}_{X/k}$) and $\underline{\Omega}^*_{X/k}$ compute Betti cohomology on affines: true by \cite[Theorem IV.1.1]{HartshorneAlgdR} for the former, and \cite{GrothendieckAlgdR} (plus cohomological descent) for the latter.
\end{proof}

\begin{remark}
	If $X$ is a smooth $k$-scheme, then the maps $a$, $b$, and $c$ from Proposition \ref{prop:derhodgefilt} are filtered isomorphisms. 
\end{remark}

Proposition \ref{prop:derhodgefilt} yields the following bound on the derived Hodge filtration:

\begin{corollary}
	Let $X$ be a finite type $k$-scheme, let $N$ be the maximal of the embedding dimensions of $\calO_{X,x}$ over all points $x \in X$. Then $\Fil^p_H(H^*(X,\widehat{\dR}_{X/k})) = 0$ for $p > N$.
\end{corollary}
\begin{proof}
If $p > N$, then the assumption on $N$ shows that $\Fil^p(\Omega^*_{X/k}) = 0$ where $\Fil$ refers to the stupid filtration. Hence,  the map $\widehat{\dR}_{X/k} \to \Omega^*_{X/k}$ from Proposition \ref{prop:derhodgefilt} factors as 
	\[ \widehat{\dR}_{X/k} \stackrel{f}{\to} \widehat{\dR}_{X/k}/\Fil^p_H \stackrel{g}{\to} \Omega^*_{X/k}.\]
	Proposition \ref{prop:derhodgefilt} shows that $g \circ f$ is a direct summand. In particular, $H^*(f)$ is injective. Since $\Fil^p_H(H^*(X,\widehat{\dR}_{X/k}))$ is exactly the kernel of $H^*(f)$, the claim follows.
\end{proof}

For a $\C$-variety $X$, the existence of Chern classes in derived de Rham cohomology \cite[\S V]{IllusieCC1} shows that the image of $c_p:K^0(X) \to H^{2p}(X^\an,\C)$ lies inside $\Fil^p_H(H^{2p}(X^\an,\C)) \cap H^{2p}(X^\an,\Q)$ where $\Fil^*_H$ is the derived Hodge filtration; here we implicitly use Corollary \ref{cor:ddrgeomcomp} and \cite[Theorem IV.1.1]{HartshorneAlgdR} to identify derived de Rham cohomology with Betti cohomology.  Proposition \ref{prop:derhodgefilt} shows that $\Fil^i_H \subset \Fil^i_\inf \subset \Fil^i_{HD}$ for all $i$. By \cite[Example 4.5]{ArapuraKangdR}, the induced inclusion 
\[ \Fil^p_\inf(H^{2p}(X^\an,\C)) \subset \Fil^p_{HD}(H^{2p}(X^\an,\C))  \]
can be strict, i.e., the infinitesimal Hodge filtration imposes more constraints on the position of Chern classes than the Hodge-Deligne filtration. To see if derived geometry does more, we ask:

\begin{question}
	Find an example of a proper $\C$-variety $X$ where the inclusion
	\[ \Fil^p_H(H^{2p}(X^\an,\C)) \cap H^{2p}(X^\an,\Q) \subset \Fil^p_\inf(H^{2p}(X^\an,\C)) \cap H^{2p}(X^\an,\Q)\] 
	is strict.
\end{question}
%\rmk{To do: find examples where all the above inclusions of filtrations are strict; study the situation for hypersurfaces.}

The preceding discussion is motivated by the Hodge conjecture which predicts that, on a smooth proper $\C$-variety, $\Q$-cohomology classes of degree $2p$ that lie in the $p^{\mathrm{th}}$-piece of the Hodge-Deligne filtration come from algebraic cycles. The same statement for singular varieties is false: Bloch gave a counterexample in \cite[Appendix A]{JannsenMixedMotives} using Mumford's work \cite{Mumfordzerocycle} on $0$-cycles on surfaces. In view of Proposition \ref{prop:derhodgefilt} and Chern classes in derived de Rham cohomology, a slightly less naive formulation\footnote{The book \cite{JannsenMixedMotives} does formulate a generalised Hodge conjecture for arbitrary varieties with no known counterexamples. However, this conjecture is formulated in terms of homology and is equivalent to the generalised Hodge conjecture for smooth projective varieties. In particular, to the best of our knowledge, it does not describe the algebraic classes in $H^{2p}(X^\an,\C)$ for a singular proper variety $X$.} of the Hodge conjecture for singular varieties is obtained by replacing the Hodge-Deligne filtration with the derived Hodge filtration. Unfortunately, this formulation is also wrong: we show below that Bloch's example also works in this setting. We do not know any counterexamples defined over $\overline{\Q}$.

%This computation also negatively answers the question raised in \cite[Remark 3.13]{ArapuraKangdR}. 

\begin{example}
	\label{ex:blochsingularhodge}
	Let $S_0 \subset \P^3_{\C}$ be a smooth degree $d \geq 4$ surface that is defined over $\overline{\Q}$. Let $p \in S_0(\C)$ be a point that is generic over $\overline{\Q}$. Set $S = \Bl_p(S_0)$, $\P = \Bl_p(\P^3)$, and let $X$ be obtained by glueing $\P$ with itself transversally along $S$. Then we will show:
	\begin{enumerate}
		\item $\Fil^2_H(H^4(X^\an,\C)) = H^4(X^\an,\C)$
		\item $c_2:K^0(X) \otimes {\Q} \to \Fil^2_H(H^4(X^\an,\C)) \cap H^4(X^\an,\Q)$ is not surjective. 
	\end{enumerate}
	Granting (1), the target of the map $c_2$ above is simply $H^4(X^\an,\Q)$, so (2) follows from Bloch's observation that $\dim(c_2(K^0(X) \otimes \Q)) = 2$ (via the genericity of $p$, see \cite[Appendix A]{JannsenMixedMotives}), while $\dim(H^4(X^\an,\Q)) = 3$ (via a Mayer-Vietoris sequence). To show (1), we must show that the map 
	\[ H^4(X,\widehat{\dR}_{X/\C}) \to H^4(X,\calO_X \stackrel{d}{\to} L_{X/\C}) \]
	is $0$. The Mayer-Vietoris sequence shows $H^4(X,\calO_X) = 0$, so it suffices to check that $H^3(X,L_{X/\C}) = 0$. Let $i_1,i_2:\P \to X$ be the two given inclusions, and let $i:S \to X$ be the composite of either $i_1$ or $i_2$ with the inclusion $S \subset \P$. Then $L_{X/\C} \simeq \Omega^1_{X/\C}$ as $X$ is lci and reduced. To show $H^3(X,\Omega^1_{X/\C}) = 0$, note that there is a sequence 
	\[ \Omega^1_{X/\C} \stackrel{a}{\to} i_{1,*} \Omega^1_{\P/\C} \oplus i_{2,*} \Omega^1_{\P/\C} \stackrel{b}{\to} i_* \Omega^1_{S/\C}\]
	with $a$ being the canonical map, and $b$ being the difference of the two canonical maps. Clearly $b$ is surjective. A check in local co-ordinates shows\footnote{The claim $\im(a) = \ker(b)$ is \'etale local on $X$, and makes sense for any ``normal crossings union'' of smooth varieties. Thus,  we may assume that $X$ is the product of a nodal curve with a smooth variety. The claim is easy to check for nodal curves, and then follows for $X$ by the Kunneth formula for $\Omega^1_{-/\C}$ as well as the exact sequence $1 \to \calO_X \to \calO_{Z_1} \oplus \calO_{Z_2} \to \calO_{Z_1 \cap Z_2} \to 1$.} that $\im(a) = \ker(b)$. This gives an exact sequence
	\[ 1 \to \ker(a) \to \Omega^1_{X/\C} \to \ker(b) \to 1.\]
It now suffices to show $H^3(X,\ker(a)) = H^3(X,\ker(b)) = 0$.  Since $a$ is an isomorphism away from $S$ and $\dim(S) = 2$, it is clear that $H^3(X,\ker(a)) = 0$. For $\ker(b)$, since $H^3(\P,\Omega^1_{\P/\C}) = 0$, it suffices to check that $H^2(S,\Omega^1_{S/\C}) = 0$ or, equivalently, that $H^1(S,\omega_S) =  0$ (by Hodge symmetry). If $\pi:S \to S_0$ denotes the blowup, then $\R\pi_* \omega_S \simeq \omega_{S_0}$ as $S_0$ has rational singularities, so we want $H^1(S_0,\omega_{S_0}) = 0$. The adjunction identification $\omega_{S_0} \simeq \omega_{\P^3}|_{S_0} \otimes (I_{S_0}/I^2_{S_0})^\vee \simeq \calO(d-4)|_{S_0}$ does the rest as $H^i(\P^3,\calO(j)) = 0$ for $i \in \{1,2\}$ and any $j$.
\end{example}

\begin{remark}
In Example \ref{ex:blochsingularhodge}, we showed that $H^4(X,\Omega^*_{X/\C}) = \Fil^2 H^4(X,\Omega^*_{X/\C})$ where $\Fil$ refers to the stupid filtration on $\Omega^*_{X/\C}$. In particular, this example negatively answers (as predicted) the question in \cite[Remark 3.13]{ArapuraKangdR}.
\end{remark}

\begin{remark}
In Example \ref{ex:blochsingularhodge}, the Mayer-Vietoris sequence shows that $H^4(X^\an,\Q) \simeq \Q(-2)^3$ as mixed Hodge structures, so the problem encountered in that example  does not come from a mismatch of weights.
\end{remark}

\begin{remark}
	The $3$-fold $X$ in Example \ref{ex:blochsingularhodge} is the fibre over $p$ of proper flat family $\pi:\calX \to S_0$ defined over $\overline{\Q}$. The fibre over any $\overline{\Q}$-point $s$ of $S_0$ is a $3$-folds $X_s$ homotopy-equivalent to $X$  with $\dim(c_2(K^0(X_s) \otimes \Q)) = 3$ provided we assume the Bloch-Beilinson conjecture that zero cycles of degree $0$ on $S_0$ defined over $\overline{\Q}$ are torsion (since the Albanese of $S_0$ is trivial). This shows the necessity of the genericity assumption on $p$ in Example \ref{ex:blochsingularhodge}.
\end{remark}

\bibliography{my}

\end{document}